\documentclass[12pt]{article}

\usepackage[T1]{fontenc}

\usepackage{mathptmx}
\usepackage{helvet}
\usepackage{courier}

\usepackage{amsmath}
\usepackage{amsthm}
\usepackage{amsfonts}
\usepackage{amssymb}

\theoremstyle{plain}
\newtheorem{theorem}{Theorem}

\theoremstyle{definition}

\usepackage{natbib}
\setcitestyle{round}

\begin{document}

\author{Anatoliy Malyarenko\thanks{M\"{a}lardalen University, V\"{a}ster{\aa}s, Sweden}}

\title{Spectral expansions of cosmological fields}

\date{\today}

\maketitle

\begin{abstract}
We give a review of the theory of random fields defined on the observable
part of the Universe that satisfy the cosmological principle, i.e.,
invariant with respect to the $6$-dimensional group $\mathcal{G}$ of the
isometries of the time slice of the
Friedmann--Lema\^{\i}tre--Robertson--Walker standard chart. Our new results
include proof of spectral expansions of scalar and spin weighted
$\mathcal{G}$-invariant cosmological fields in open, flat, and closed
cosmological models.
\end{abstract}

\section{Introduction}

The space mission Euclid, planned for launch in 2019, is expected, just to
mention a few, to measure the shapes of 1.5 billion galaxies. In mathematical
terms, Euclid will observe a single realisation of a random section of a spin
bundle over a ball of radius about 10 billion light years. Analysis of Euclid
data should base on spectral expansions of different cosmological fields. In
this paper, we make attempt to review the known expansions and to present
some new results.

In Section~\ref{sec:preliminaries}, we review the essential axioms,
definitions, and results of the Standard Cosmological Model. We consider
open, flat, and closed Friedmann--Lema\^{\i}tre--Robertson--Walker
cosmological models. The time slice $X$ of each of the above models has
$6$-dimensional group of isometries $\mathcal{G}$, where
$\mathcal{G}=\mathrm{SO}_0(1,3)$ for the open model,
$\mathcal{G}=\mathrm{ISO}(3)$ for the flat model, and
$\mathcal{G}=\mathrm{SO}(4)$ for the closed model.

Section~\ref{sec:scalar} introduces scalar cosmological fields, both
deterministic and stochastic. In particular,
Subsection~\ref{sub:deterministic} reviews expansions of deterministic
cosmological fields in eigenfunction of the $\mathcal{G}$-invariant
Laplace--Beltrami operator. The radial coefficients $R_{k\ell}(\chi)$ and
$R_{\omega\ell}(\chi)$ of the above expansions are used in
Subsection~\ref{sub:stochastic} to deduce spectral expansions of
$\mathcal{G}$-invariant stochastic fields. The above expansions are proved in
Theorem~\ref{th:1} similarly for all three models, using tools from harmonic
analysis.

In Section~\ref{sec:spin} we first motivate introduction of spin stochastic
cosmological fields, using the theory of the effective lensing potential.
After short review of spin weighted spherical harmonics, we prove
Theorem~\ref{th:2}, that gives the spectral expansion of the convergence,
shear, and of the two gravitational flexions cosmological fields.

\section{Preliminaries from cosmology}\label{sec:preliminaries}

According to modern cosmology, the Universe is governed by Einstein's
equations. Solution of the above equations is beyond the capabilities of
modern computers. Cosmologists use another approach. They start from a
simple exactly solvable model as a background solution and then add
perturbations order by order.

The background solution is described as follows. The Universe is a
$4$-di\-men\-si\-o\-nal Lorentzian manifold $M$ satisfying the following two
axioms.

\begin{enumerate}

\item \textbf{The cosmological principle}: at each epoch, the Universe
    presents the same aspect from every point.

\item \textbf{The Weyl's postulate} after \cite{Weyl1923}:

\begin{quote}
``The particles of the substratum (representing the nebulae) lie, in
spacetime of the cosmos, on a bundle of geodesics diverging from a point
in the (finite or infinite) past.''
\end{quote}

\end{enumerate}

The ``substratum'' is the underlying fluid defining the overall kinematics of
a system of galaxies. The galaxies follow their geodesics. Only one geodesic
is passing through each point in $M$, except at the origin. This allows one
to define the concept of \emph{fundamental observer}, one for each geodesic.
Each of these is carrying a standard clock, for which they can synchronise
and fix a \emph{cosmic time}.

For each fundamental observer, there exists a \emph{standard chart}
$(t,\chi,\theta,\varphi)$ of the manifold~$M$ centred at the observer in
which the metric has the \emph{Friedmann--Lema\^{\i}tre--Robertson--Walker}
form:
\begin{equation}\label{eq:metric}
\mathrm{d}s^2=-c^2\,\mathrm{d}t^2+a^2(t)[\mathrm{d}\chi^2+f^2_K(\chi)
(\mathrm{d}\theta^2+\sin^2\theta\,\mathrm{d}\varphi^2)],
\end{equation}
where the function $f_K$ depends on the curvature $K$ as
\[
f_K(\chi)=
\begin{cases}
(-K)^{-1/2}\sinh(\sqrt{-K}\chi),&K<0,\\
\chi,&K=0,\\
K^{-1/2}\sin(\sqrt{K}\chi),&K>0,
\end{cases}
\]
and where $t$ is the cosmic time, $\chi$ the \emph{comoving distance},
$\theta$ and $\varphi$ are spherical coordinates on the celestial sphere, $c$
is the speed of light in vacuum, $a(t)$ is the \emph{scale factor}. The
comoving distance from the observer to a galaxy is the distance the galaxy
will have from the observer when it is as old as the observer now. The
physical distance is the product of the scale factor and the comoving
distance. We choose the scale factor to be dimensionless and set $a(t_0)=1$,
where $t_0$ is the value of the cosmic time now. Thus, the comoving distance
$\chi$ has the dimension of length, while $K$ has the dimension of inverse
length squared. The velocity of a galaxy with respect to a fundamental
observer is called the \emph{peculiar velocity}.

The time slice $D_t$ of the domain $D$ of the above chart is the Riemannian
manifold $X$ with metric
$a^2(t)[\mathrm{d}\chi^2+f^2_K(\chi)(\mathrm{d}\theta^2+\sin^2\theta\,\mathrm{d}\varphi^2)]$.
The geodesic lines of fundamental observers are orthogonal to the time
slices. The cosmological principle implies that each time slice $X$ is the
$3$-dimensional space of constant curvature. When the curvature is negative,
the model is called \emph{open}, the group of isometries of $X$ is
$\mathcal{G}=\mathrm{SO}_0(1,3)$, and $X=\mathbb{H}^3$, the hyperbolic space.
When the curvature is equal to $0$, the model is called \emph{flat}, the
group of isometries of $X$ is $\mathcal{G}=\mathrm{ISO}(3)$, and
$X=\mathbb{R}^3$. Finally, when the curvature is positive, the model is
called \emph{closed}, the group of isometries of $X$ is
$\mathcal{G}=\mathrm{SO}(4)$, and $X=\mathbb{S}^3$, the sphere. For all three
models, the stationary subgroup of the point $(t,0,0,0)$ is
$\mathcal{K}=\mathrm{SO}(3)$.

The surface $S(\chi)$ of a comoving sphere of radius $\chi$ is $S(\chi)=4\pi
f^2_K(\chi)$. The $\mathcal{G}$-invariant measure on $X$ in the above chart
has the form
\[
\mathrm{d}x=f^2_K(\chi)\sin\theta\,\mathrm{d}\chi\,\mathrm{d}\theta\,\mathrm{d}\varphi.
\]
The comoving distance $\chi$ runs from $0$ to $\pi/\sqrt{K}$ when $K>0$, and
from $0$ to $\infty$ otherwise.

The scale factor, $a(t)$, is not observable. However, consider the
electromagnetic wave of length $\lambda_1$ emitted by a distant galaxy at
time $t_1<t_0$ that approaches the Earth now. It has another wave length,
$\lambda_0$. The ratio $\lambda_0/\lambda_1$ is observable. It is denoted by
\[
\frac{\lambda_0}{\lambda_1}=1+z.
\]
The quantity $z$ is the celebrated \emph{redshift}. Moreover,
\[
a(t_1)=\frac{1}{1+z}    .
\]
In 1929, Hubble discovered that the redshift is positive for all but few
close galaxies. It follows that $a(t_1)<1$. In other words, the Universe is
expanding. The galaxies without redshift are those having big enough
peculiar velocities to stand against the expansion in the standard chart.

Denote
\[
H_0=\dot{a}(t_0),
\]
where a dot denote differentiation with respect to $t$. $H_0$ is called the
\emph{Hubble constant}. The number
\[
\rho_c=\frac{3H^2_0}{8\pi G},
\]
where $G$ is the Newton's gravitational constant, is called the
\emph{critical density}.

In general, $H(z)$ denote the Hubble constant as measured by an imaginary
astronomer working at redshift $z$. The expansion of the Universe is
governed by the \emph{Friedmann equation}:
\[
H^2(z)=H^2_0[\Omega_R(1+z)^4+\Omega_M(1+z)^3+\Omega_K(1+z)^2+\Omega_{\Lambda}],
\]
where the parameters are as follows.

$\Omega_M$ is the \emph{density parameter for non-relativistic matter}
(which moves with speed much less than $c$ with respect to the fundamental
observer). We have
\[
\Omega_M=\frac{\rho_{M,t_0}}{\rho_c},
\]
where $\rho_{M,t_0}$ is the density of the non-relativistic matter now, at
$t=t_0$.

$\Omega_{\Lambda}$ is the \emph{vacuum energy parameter}:
\[
\Omega_{\Lambda}=\frac{\Lambda c^2}{3H^2_0},
\]
where $\Lambda$ is the \emph{cosmological constant}.

$\Omega_K$ is the \emph{curvature parameter}:
\[
\Omega_K=-\frac{Kc^2}{H^2_0}.
\]

Finally, $\Omega_R$ is the \emph{radiation density parameter}. We have
\[
\Omega_R=\frac{\rho_{R,t_0}}{\rho_c},
\]
where $\rho_{R,t_0}$ is the density of the radiation now. When $z=0$, we
obtain
\[
\Omega_R+\Omega_M+\Omega_K+\Omega_{\Lambda}=1.
\]

The values for the above parameters, adapted from \cite{Efstathiou2013}, are
shown in Table~\ref{tab:1}, where $1$ Mpc = $3.09\times 10^{26}$ m, and the
quoted errors show the $68$\% confidence level. The contribution from
radiation is negligibly small, $\Omega_R=(4.9\pm 0.5)\times 10^{-5}$.

\begin{table}
\centering
\begin{tabular}{|l|l|}
\hline \textbf{Parameter} & \textbf{Value} \\
\hline $H_0$ & $67.80\pm 0.77$ km/s/Mpc \\
$\Omega_M$ & $0.315^{+0.016}_{-0.018}$ \\
$\Omega_K$ & $-0.0010^{+0.0062}_{-0.0065}$ \\
$\Omega_{\Lambda}$ & $0.685^{+0.018}_{-0.016}$ \\
\hline
\end{tabular}
\caption{Cosmological parameters}\label{tab:1}
\end{table}

The comoving distance $\chi$ from the fundamental observer at $z=0$ to an
object with redshift $z$ is
\[
\chi(z)=\frac{c}{H_0}\int^z_0\frac{\mathrm{d}u}{\sqrt{\Omega_R(1+u)^4
+\Omega_M(1+u)^3+\Omega_K(1+u)^2+\Omega_{\Lambda}}}.
\]
The \emph{look-back time} $t_L$ to the above object is the difference
between $t_0$ and the cosmic time at the moment when photons were emitted by
the object. It is calculated as
\[
t_L=\frac{1}{H_0}\int^z_0\frac{\mathrm{d}u}{(1+u)\sqrt{\Omega_R(1+u)^4
+\Omega_M(1+u)^3+\Omega_K(1+u)^2+\Omega_{\Lambda}}}.
\]

In Friedmann equation, $a=0$ is reached after a finite time. Near $a=0$, the
Hubble parameter $H(z)$ explodes. This point is therefore called \emph{Big
Bang}.

\section{Scalar fields}\label{sec:scalar}

\subsection{Deterministic fields}\label{sub:deterministic}

Let $f\colon M\to\mathbb{R}$ be a deterministic scalar cosmological field.
On each time slice $X$, it is natural to decompose the restriction of $f$ to
$X$ in eigenfunctions of $\mathcal{G}$-invariant differential operators.

Let $d(x,y)$ be a distance between the points $x$ and $y\in X$. In all three
models, $X$ is a \emph{two-point homogeneous space}. That is: for any
two-point pairs $x$, $y\in X$, $x'$, $y'\in X$ satisfying $d(x,y)=d(x',y')$,
there exists an isometry $g\in\mathcal{G}$ such that $gx=x'$, $gy=y'$. By
\cite[Chapter~2, Proposition~4.11]{Helgason1984}, all
$\mathcal{G}$-invariant differential operators are polynomials in the
Laplace--Beltrami operator $\Delta_K$ of the Riemannian manifold $X$. In
what follows we call $\Delta_K$ just the Laplacian. We have
\[
\Delta_K=\frac{1}{f^2_K(\chi)}\frac{\partial}{\partial\chi}\left(f^2_K(\chi)
\frac{\partial}{\partial\chi}\right)+\frac{1}{f^2_K(\chi)}\Delta,\\
\]
where
\[
\Delta=\frac{1}{\theta^2}\frac{\partial}{\partial\theta}\left(\theta^2
\frac{\partial}{\partial\theta}\right)
+\frac{1}{\sin^2\theta}\frac{\partial^2}{\partial\varphi^2}
\]
is the Laplacian on the sphere $\mathbb{S}^2$.

The eigenvalue problem for the Laplacian $\Delta_K$ is traditionally written
as the \emph{Helmholtz equation}
\[
\Delta_KQ_{\mathbf{k}}(x)=-(k^2-K)Q_{\mathbf{k}}(x),
\]
where $\mathbf{k}$ is called the \emph{wave vector} and labels the
\emph{modes}, $Q_{\mathbf{k}}(x)$, and where $k$ is called the \emph{wave
number}. The wave numbers take the values in $[0,\infty)$ for $K\leq 0$ and
in the set $\{\,(\omega+1)\sqrt{K}\colon\omega\in\mathbb{Z}_+\,\}$
otherwise, where $\mathbb{Z}_+$ is the set of all nonnegative integers.
Different authors use different normalisation of modes.

The Helmholtz equation may be solved by the method of separation of
variables. We use the normalisation by \cite{Peter2009}:
\[
Q_{k\ell m}(\chi,\theta,\varphi)=(2\pi)^{3/2}R_{k\ell}(\chi)
Y_{\ell m}(\theta,\varphi),
\]
where $Y_{\ell m}(\theta,\varphi)$ are \emph{spherical harmonics}, the
eigenfunctions of $\Delta$:
\[
\Delta Y_{\ell m}=-\ell(\ell+1)Y_{\ell m},\qquad\ell\geq 0,\quad-\ell\leq m\leq\ell.
\]
Spherical harmonics are related to \emph{associated Legendre polynomials}
$P^m_{\ell}$ through
\[
Y_{\ell m}(\theta,\varphi)=\sqrt{\frac{(2\ell+1)(\ell-m)!}{4\pi(\ell+m)!}}
P^m_{\ell}(\cos\theta)e^{\mathrm{i}m\varphi}.
\]
The radial part of the Helmholtz equation becomes
\begin{equation}\label{eq:radial}
\frac{1}{f^2_K(\chi)}\frac{\mathrm{d}}{\mathrm{d}\chi}\left(f^2_K(\chi)
\frac{\mathrm{d}R_{k\ell}(\chi)}{\mathrm{d}\chi}\right)+\left[k^2-K
-\frac{\ell(\ell+1)}{f^2_K(\chi)}\right]R_{k\ell}(\chi)=0.
\end{equation}

When $K<0$, put $\omega=k/\sqrt{-K}$ and
\[
N_{k\ell}=\prod^{\ell}_{n=0}(\omega^2+n^2).
\]
The solution to \eqref{eq:radial} has the form
\begin{equation}\label{eq:negativeR}
R_{k\ell}(\chi)=\sqrt{\frac{\pi N_{k\ell}}{2\omega^2\sinh(\sqrt{-K}\chi)}}
P^{-1/2-\ell}_{-1/2+\mathrm{i}\omega}(\cosh(\sqrt{-K}\chi)),
\end{equation}
where $P^{-1/2-\ell}_{-1/2+\mathrm{i}\omega}$ is the \emph{associated
Legendre function of the first kind}.

When $K=0$, the solution to \eqref{eq:radial} has the form
\[
R_{k\ell}(\chi)=\sqrt{\frac{2}{\pi}}j_{\ell}(k\chi),
\]
where $j_{\ell}$ is the \emph{spherical Bessel function}. The dimensionless
distance in this case is $r=k\chi$.

Finally, when $K>0$, put $\omega=k/\sqrt{K}-1$ and
\[
M_{\omega\ell}=\prod^{\ell}_{n=0}((\omega+1)^2-n^2).
\]
Note that $M_{\omega\ell}=0$ for $\ell>\omega$. The solution to
\eqref{eq:radial} has the form
\begin{equation}\label{eq:positiveR}
R_{\omega\ell}(\chi)=\sqrt{\frac{\pi M_{\omega\ell}}{2(\omega+1)^2\sin(\sqrt{K}\chi)}}
P^{-1/2-\ell}_{1/2+\omega}(\cos(\sqrt{K}\chi)),
\end{equation}
see \cite[Equation~(A21)]{Abbott1986}.

For $K\leq 0$, the expansion of a function $f\in L^2(X,\mathrm{d}x)$ into
modes has the form
\begin{equation}\label{eq:deterministic}
f(\chi,\theta,\varphi)=(2\pi)^{-3/2}\sum^{\infty}_{\ell=0}\sum^{\ell}_{m=-\ell}
\int^{\infty}_0k^2f_{\ell m}(k)Q_{k\ell m}(\chi,\theta,\varphi)\,\mathrm{d}k.
\end{equation}
For $K>0$, the above expansion reads
\[
f(\chi,\theta,\varphi)=\left(\frac{K}{2\pi}\right)^{3/2}\sum^{\infty}_{\omega=0}
\sum^{\omega}_{\ell=0}\sum^{\ell}_{m=-\ell}(\omega+1)^2f_{\ell m}(k)
Q_{k\ell m}(\chi,\theta,\varphi),
\]
where the Fourier coefficients, $f_{\ell m}(k)$, are given by
\[
f_{\ell m}(k)=(2\pi)^{-3/2}\int_Xf(\chi,\theta,\varphi)\overline{Q_{k\ell m}
(\chi,\theta,\varphi)}f^2_K(\chi)\sin\theta\,\mathrm{d}\chi\,\mathrm{d}\theta
\,\mathrm{d}\varphi.
\]

Assume that $K<0$ and $f$ depends only on $\chi$. Then we have
\[
Q_{k00}(\chi)=\pi^{3/2}\frac{1}{\sqrt{\sin(\sqrt{-K}\chi)}}
P^{-1/2}_{-1/2+\mathrm{i}\omega}(\cosh(\sqrt{-K}\chi)).
\]
Using the formula
\[
P^{-1/2}_{-1/2+\mathrm{i}\omega}(\cosh r)=\sqrt{\frac{2}{\pi\sinh r}}\frac{\sin(\omega r)}{\omega}
\]
we obtain
\[
Q_{k00}(\chi)=\pi\sqrt{2}\frac{\sqrt{-K}\sin(k\chi)}{k\sinh(\sqrt{-K}\chi)},
\]
and \eqref{eq:deterministic} becomes
\begin{subequations}\label{eq:negative}
\begin{equation}\label{eq:negative:a}
f(\chi)=\frac{1}{2\sqrt{\pi}}\int^{\infty}_0f_{00}(k)\frac{\sqrt{-K}\sin(k\chi)}
{k\sinh(\sqrt{-K}\chi)}k^2\,\mathrm{d}k,
\end{equation}
where
\begin{equation}\label{eq:negative:b}
f_{00}(k)=\frac{1}{2\sqrt{\pi}}\int^{\infty}_0f(\chi)\frac{\sqrt{-K}\sin(k\chi)}
{k\sinh(\sqrt{-K}\chi)}S(\chi)\,\mathrm{d}\chi.
\end{equation}
\end{subequations}

Similar calculations show that in the case of $K=0$ we have
\begin{subequations}\label{eq:zero}
\begin{equation}\label{eq:zero:a}
f(\chi)=\frac{1}{\pi\sqrt{2}}\int^{\infty}_0f_{00}(k)\frac{\sin(k\chi)}{k\chi}
k^2\,\mathrm{d}k,
\end{equation}
where
\begin{equation}\label{eq:zero:b}
f_{00}(k)=\frac{1}{\pi\sqrt{2}}\int^{\infty}_0f(\chi)\frac{\sin(k\chi)}{k\chi}
S(\chi)\,\mathrm{d}\chi,
\end{equation}
\end{subequations}
while in the case of $K>0$ we have
\begin{subequations}\label{eq:positive}
\begin{equation}\label{eq:positive:a}
f(\chi)=\frac{K^{3/2}}{2\sqrt{\pi}}\sum^{\infty}_{\omega=0}f_{00}(\omega)
\frac{\sin((\omega+1)\sqrt{K}\chi)}{(\omega+1)\sin(\sqrt{K}\chi)}(\omega+1)^2,
\end{equation}
where
\begin{equation}\label{eq:positive:b}
f_{00}(\omega)=\frac{1}{2\sqrt{\pi}}\int^{\pi/\sqrt{K}}_0f(\chi)
\frac{\sin((\omega+1)\sqrt{K}\chi)}{(\omega+1)\sin(\sqrt{K}\chi)}S(\chi)\,\mathrm{d}\chi.
\end{equation}
\end{subequations}

The above expansions, with various normalisations, were considered by
\cite{Fock1935}, \cite{Schroedinger1939,Schroedinger1957},
\cite{Lifshitz1963}, \cite{Harrison1967}, \cite{Abbott1986}, among others.

\subsection{Stochastic fields}\label{sub:stochastic}

Let $f\colon M\times\Omega\to\mathbb{C}$ be a stochastic scalar cosmological
field defined on a probability space $(\Omega,\mathfrak{F},\mathsf{P})$. We
assume that $f$ has a finite variance, $\mathsf{E}[|f(t,x)|^2]<\infty$ and
that $f$ is mean-square continuous, i.e., the map $M\to L^2(\Omega)$,
$(t,x)\mapsto f(t,x)$ is continuous, where $L^2(\Omega)$ is the Hilbert
space of all random variables with finite variance.

Consider the restriction $f(x)$ of the field $f(t,x)$ to the time slice $X$.
By the cosmological principle, the field $f(x)$ must be
$\mathcal{G}$-invariant, i.e., for any positive integer $n$, for any
distinct points $x_1$, \dots, $x_n\in X$, and for any $g\in\mathcal{G}$,
random vectors $(f(x_1),\dots,f(x_n))^{\top}$ and
$(f(gx_1),\dots,f(gx_n))^{\top}$ must have the same distribution.

In particular, the expected value $\mathsf{E}[f(x)]$ is a constant, say
$a\in\mathbb{C}$, and the autocorrelation function
\[
R(x,y)=\mathsf{E}[(f(x)-a)\overline{(f(y)-a)}]
\]
is $\mathcal{G}$-invariant, i.e., $R(gx,gy)=R(x,y)$. Such fields are called
\emph{wide-sense homogeneous}. Note that $R(x,y)$ is a positive-definite
function.

The spectral theory of wide-sense homogeneous random fields on topological
groups and homogeneous spaces has been developed by \cite{yaglom:61}. He
proved spectral expansions of wide-sense homogeneous random fields on
separable topological groups of type I and their homogeneous spaces. A
topological group of type I is such a group $\mathcal{G}$ that every unitary
representation of $\mathcal{G}$ generates an operator algebra of type I
(see, for example, \cite{Naimark1972}). All three groups
$\mathrm{SO}_0(1,3)$, $\mathrm{ISO}(3)$, and $\mathrm{SO}(4)$, are of type
I.

Moreover, the above groups share the following property. Let $U$ be an
irreducible unitary representation of the group $\mathcal{G}$. Consider the
restriction of $U$ to the compact subgroup $\mathcal{K}=\mathrm{SO}(3)$. The
above restriction is equivalent to a direct sum of irreducible unitary
representations of $\mathcal{K}$. The multiplicity of the trivial
representation of the group $\mathcal{K}$ in the above sum is equal either
to $0$ or to $1$.

Let $\hat{\mathcal{G}}_{\mathcal{K}}$ be the set of all equivalence classes
of irreducible unitary representations of the group $\mathcal{G}$ for which
the above multiplicity is equal to $1$. Let
$U^{\omega}\in\hat{\mathcal{G}}_{\mathcal{K}}$, and let $H_{00}$ be the
one-dimensional complex Hilbert space in which the trivial representation of
the group $\mathcal{K}$ is realised. Let $\mathbf{e}$ be a vector of unit
length in $H_{00}$. It is easy to see that the function
\[
\Phi_{\omega}(g)=(U^{\omega}\mathbf{e},\mathbf{e})_{H_{00}}
\]
does not depend on the choice of the vector $\mathbf{e}$. It is called the
\emph{zonal spherical function} of the group $\mathcal{G}$.

Let $\pi\colon\mathcal{G}\to X$ be the natural projection, let $x$, $y\in
X$, and let $g_1\in\pi^{-1}(x)$, $g_2\in\pi^{-1}(y)$. By
\cite[Theorem~6']{yaglom:61}, formula
\[
R(x,y)=\int_{\hat{\mathcal{G}}_{\mathcal{K}}}\Phi_{\omega}(g^{-1}_2g_1)\,
\mathrm{d}\nu(\omega)
\]
determines a one-to-one correspondence between the autocorrelation functions
$R(x,y)$ of wide-sense homogeneous random fields on $X$ and finite measures
$\nu$ on $\hat{\mathcal{G}}_{\mathcal{K}}$. When $X$ is a two-point
homogeneous space, the zonal spherical function $\Phi_{\omega}(g^{-1}_2g_1)$
depends only on $r$, the distance between $x$ and $y$.

The set $\hat{\mathcal{G}}_{\mathcal{K}}$ and the zonal spherical function
for the case of $X=\mathbb{H}^3$ were determined by \cite{Gelfand1946}. There
exist different parametrisations of the above set. We will use the
parametrisation by \cite{kostant}. The set $\hat{\mathcal{G}}_{\mathcal{K}}$
includes the \emph{principal series} parameterised by $\omega\in[0,+\infty)$,
the \emph{supplementary series} parameterised by $\omega\in(0,\mathrm{i})$,
and trivial representation with $\omega=\mathrm{i}$. We have
\[
\Phi_{\omega}(r)=\frac{\sin(\omega r)}{\omega\sinh r}
\]
and
\begin{equation}\label{eq:Kreuin}
R(x,y)=\int^{\infty}_\mathrm{i}\frac{\sin(\omega r)}{\omega\sinh r}\,\mathrm{d}\nu(\omega),
\end{equation}
where the integral is taken over the set $(0,\mathrm{i}]\cup[0,\infty)$.
Equation \eqref{eq:Kreuin} was obtained by \cite{Kreuin1949} for the case of
a finite-dimensional hyperbolic space $\mathbb{H}^n$. To compare
\eqref{eq:Kreuin} with \eqref{eq:negative:a}, we choose $r$ in such a way
that $\omega r=k\chi$, or $r=\sqrt{-K}\chi$. We obtain
\[
R(\chi)=\frac{1}{\sqrt{-K}}\int^{\infty}_{\mathrm{i}\sqrt{-K}}\frac{\sqrt{-K}\sin(k\chi)}
{k\sinh(\sqrt{-K}\chi)}\,\mathrm{d}\nu(k).
\]

Compare the last display with \eqref{eq:negative:a}. Apart from
normalisation, we note the following difference. The expansion of a
square-integrable function includes only the representations of the principal
series of the group $\mathcal{G}$ (in physical terms, \emph{sub-curvature
modes}, when the wavelength is less than the radius of curvature). On the
other hand, the expansion of a positive-definite function includes the
representations of \emph{both} the principal and the supplementary series
(sub- and \emph{super-curvature modes}, when the wavelength is greater than
the radius of curvature). This difference was noted by \cite{Lyth1995}.
\cite{Lyth1995a} describes this situation as follows.

\begin{quote}
Mathematicians have known for almost half a century that all modes must be
included to generate the most general \emph{homogeneous Gaussian random
field}, despite the fact that any square integrable \emph{function} can be
generated using only the sub-curvature modes. The former mathematical
object, not the latter, is the relevant one for physical applications.
\end{quote}

The set $\hat{\mathcal{G}}_{\mathcal{K}}$ and the zonal spherical functions
for the case of $X=\mathbb{R}^3$ were determined by \cite{Kreuin1949}. We
have $\hat{\mathcal{G}}_{\mathcal{K}}=[0,\infty)$,
\[
\Phi_{\omega}(r)=\frac{\sin(\omega r)}{\omega r}.
\]
To compare this with \eqref{eq:zero:a}, choose $\omega=k$ and $r=\chi$. We
obtain
\[
R(\chi)=\int^{\infty}_0\frac{\sin(k\chi)}{k\chi}\,\mathrm{d}\nu(k).
\]

\cite{Cartan1929} determined the set $\hat{\mathcal{G}}_{\mathcal{K}}$ and
the zonal spherical functions for the case of $X=\mathbb{S}^3$. In this case
$\hat{\mathcal{G}}_{\mathcal{K}}=\mathbb{Z}_+$,
\[
\Phi_{\omega}(r)=\frac{\sin((\omega+1)r)}{(\omega+1)\sin r}.
\]
To compare this with \eqref{eq:positive:a}, put $r=\sqrt{K}\chi$ and recall
that $\omega=k/\sqrt{K}-1$. We obtain
\[
R(\chi)=\sum^{\infty}_{\omega=0}\nu(\omega)\frac{\sin((\omega+1)\sqrt{K}\chi)}
{(\omega+1)\sin(\sqrt{K}\chi)},
\]
where $\nu(\omega)\geq 0$ and $\sum^{\infty}_{\omega}\nu(\omega)<\infty$.

Equations \eqref{eq:negative}, \eqref{eq:zero}, and \eqref{eq:positive} have
the following group-theoretical interpretation. Consider the set of double
cosets $\mathcal{K}\backslash\mathcal{G}/\mathcal{K}$. When $K\leq 0$, this
set is the interval $[0,\infty)$, otherwise this set is $[0,\pi\sqrt{K}]$.
The measure $S(\chi)\,\mathrm{d}\chi$ is a unique up to a constant multiplier
$\mathcal{K}$-bi-invariant measure on
$\mathcal{K}\backslash\mathcal{G}/\mathcal{K}$. For any choice of the above
multiplier, say $C$, there exists a unique measure on
$\hat{\mathcal{G}}_{\mathcal{K}}$ called the \emph{Plancherel measure}. The
Plancherel measure is characterised by the following property: equations
\eqref{eq:negative}--\eqref{eq:positive} determine an isometric isomorphism
between the space of $\mathcal{K}$-bi-invariant square-integrable functions
on $\mathcal{G}$ with measure $CS(\chi)\,\mathrm{d}\chi$ and the space of
square-integrable functions on $\hat{\mathcal{G}}_{\mathcal{K}}$ with
Plancherel measure. When $K\leq 0$, the measure $k^2\,\mathrm{d}k$ is
proportional to the Plancherel measure on $\hat{\mathcal{G}}_{\mathcal{K}}$.
Note that the Plancherel measure is equal to zero for the supplementary
series $\omega\in(0,\mathrm{i})$ and the trivial representation
$\omega=\mathrm{i}$. When $K>0$, the Plancherel measure of a point
$\omega\in\hat{\mathcal{G}}_{\mathcal{K}}$ is proportional to $(\omega+1)^2$.
Equations \eqref{eq:negative:b}, \eqref{eq:zero:b}, and \eqref{eq:positive:b}
define the \emph{spherical Fourier transform}, while equations
\eqref{eq:negative:a}, \eqref{eq:zero:a}, and \eqref{eq:positive:a} define
the inverse transform.

To find the spectral expansion of the random field $f$, we have to find a
measurable space $(V,\mathfrak{V},\mu)$ and a function $h\colon X\times
V\to\mathbb{C}$ such that
\begin{equation}\label{eq:Kahrounen}
R(x,y)=\int_Vh(x,v)\overline{h(y,v)}\,\mathrm{d}\mu(v).
\end{equation}
In this case there exists a complex-valued scattered random measure $Z$ on
$V$ with $\mu$ as its control measure, i.e., for any $A$, $B\in\mathfrak{V}$
we have
\[
\mathsf{E}[Z(A)\overline{Z(B)}]=\mu(A\cap B)
\]
such that
\[
f(x)=\int_Vh(x,v)\,\mathrm{d}Z(v).
\]
Note that if $h(x,v)$ satisfies \eqref{eq:Kahrounen}, then
$h(x,v)e^{\mathrm{i}\varphi(v)}$ also satisfies the above equation.

For any $\omega\in\hat{\mathcal{G}}_{\mathcal{K}}$, consider the restriction
of the representation $U^{\omega}$ to $\mathcal{K}$. In the cases of the
open model with $\omega\neq\mathrm{i}$ and the flat model with $\omega\neq
0$, the above restriction is equivalent to the direct sum of the irreducible
unitary representations $D^{\ell}$ of the group $\mathcal{K}$ for all
$\ell\in\mathbb{Z}_+$ acting in Hilbert spaces $H_{\ell}$. In the case of
the closed model, the above restriction is equivalent to the direct sum of
the representations $D^{\ell}$, $0\leq\ell\leq\omega$. The restriction of
any representation $D^{\ell}$ to the subgroup
$\mathrm{SO}(2)\subset\mathcal{K}$ is equivalent to the direct sum of the
representations $\varphi\mapsto e^{\mathrm{i}m\varphi}$ of the group
$\mathrm{SO}(2)$, $-\ell\leq m\leq\ell$. We represented the space
$H_{\omega}$ of the representation $U^{\omega}$ as the direct sum of
one-dimensional subspaces $H_{\ell m}$. Choose a vector $\mathbf{e}_{\ell
m}\in H_{\ell m}$ of unit length in such a way, that the entries of matrices
\[
D^{\ell}_{mn}=(D^{\ell}\mathbf{e}_{\ell m},\mathbf{e}_{\ell n})_{H_{\ell}}
\]
are Wigner $D$-functions:
\begin{equation}\label{eq:Wigner}
\begin{aligned}
D^{\ell}_{mn}(\varphi,\theta,\psi)&=e^{-\mathrm{i}(m\varphi+n\psi)}(-1)^m\sqrt{\frac{(\ell+m)!
(\ell-m)!}{(\ell+n)!(\ell-n)!}}\sin^{2\ell}(\theta/2)\\
&\quad\times\sum^{\min\{\ell+m,\ell+n\}}_{s=\max\{0,m+n\}}\binom{\ell+n}{s}
\binom{\ell-n}{s-m-n}(-1)^{\ell-s+n}\cot^{2s-m-n}(\theta/2).
\end{aligned}
\end{equation}

Let $g=k_1rk_2$ be the standard polar decomposition of the elements of the
group~$\mathcal{G}$. Choose $g_i=k_ir_i\in\pi^{-1}(x_i)$, $i=1$, $2$. Then
we have
\[
R(x,y)=\int_{\hat{\mathcal{G}}_{\mathcal{K}}}\Phi_{\omega}(r^{-1}_2k^{-1}_2k_1r_1)\,
\mathrm{d}\nu(\omega).
\]
By matrix multiplication
\[
R(x,y)=\sum_{\ell}\sum^{\ell}_{m,n=-\ell}D^{\ell}_{mn}(k^{-1}_2k_1)
\int_{\hat{\mathcal{G}}_{\mathcal{K}}}U^{\omega}_{\ell,m;0,0}(r_1)
\overline{U^{\omega}_{\ell,n;0,0}(r_2)}\,\mathrm{d}\nu(\omega).
\]
The matrix entry $U^{\omega}_{\ell,m;0,0}(r_1)$ is called the
\emph{associated spherical function}. For all three models, they are
different from $0$ only when $m=0$. Therefore,
\[
R(x,y)=\sum_{\ell}D^{\ell}_{00}(k^{-1}_2k_1)\int_{\hat{\mathcal{G}}_{\mathcal{K}}}
U^{\omega}_{\ell,0;0,0}(r_1)\overline{U^{\omega}_{\ell,0;0,0}(r_2)}\,\mathrm{d}\nu(\omega).
\]
Again by matrix multiplication
\[
R(x,y)=\sum_{\ell}\sum^{\ell}_{m=-\ell}D^{\ell}_{m0}(k_1)\overline{D^{\ell}_{m0}(k_2)}
\int_{\hat{\mathcal{G}}_{\mathcal{K}}}U^{\omega}_{\ell,0;0,0}(r_1)
\overline{U^{\omega}_{\ell,0;0,0}(r_2)}\,\mathrm{d}\nu(\omega).
\]
Wigner $D$-functions are related to spherical harmonics by
\begin{equation}\label{eq:spherical}
Y_{\ell m}(\theta,\varphi)=\sqrt{\frac{2\ell+1}{4\pi}}D^{(\ell)}_{-m0}(\varphi,\theta,0).
\end{equation}
Therefore we have
\begin{equation}\label{eq:correlation}
\begin{aligned}
R(x,y)&=4\pi\sum_{\ell}\frac{1}{2\ell+1}\sum^{\ell}_{m=-\ell}
Y_{\ell m}(\theta_1,\varphi_1)\overline{Y_{\ell m}(\theta_2,\varphi_2)}\\
&\quad\times\int_{\hat{\mathcal{G}}_{\mathcal{K}}}U^{\omega}_{\ell,0;0,0}(r_1)
\overline{U^{\omega}_{\ell,0;0,0}(r_2)}\,\mathrm{d}\nu(\omega),
\end{aligned}
\end{equation}
where
\[
r_i=
\begin{cases}
\sqrt{|K|}\chi_i,&K\neq 0,\\
k\chi_i,&K=0,
\end{cases}
\]
and $(\chi_1,\theta_1,\varphi_1)$ (resp. $(\chi_2,\theta_2,\varphi_2)$) are
coordinates of the point $x$ (resp. $y$) in the standard chart.

In the case of the open model we have
\[
V=\{\,((0,\mathrm{i}]\cup[0,\infty),\ell,m)\colon\ell\in\mathbb{Z}_+,
-\ell\leq m\leq\ell\,\},
\]
i.e., the union of countably many copies of the set
$(0,\mathrm{i}]\cup[0,\infty)$ indexed by pairs $(\ell,m)$, $\mathfrak{V}$
is the $\sigma$-field of Borel sets in $V$, and the restriction of the
measure $\mu$ to any copy is equal to $\nu$.

The associated spherical functions for this case were calculated by
\cite{Vilenkin1958}:
\[
U^{\omega}_{\ell,0;0,0}(r)=\frac{(-1)^{\ell}\sqrt{(2\ell+1)\pi}\Gamma(\mathrm{i}\omega)}
{\sqrt{2}\Gamma(\mathrm{i}\omega-\ell)}\sinh^{-1/2}r P^{-1/2-\ell}_{-1/2+\mathrm{i}\omega}(\cosh r).
\]
Substitute this formula to \eqref{eq:correlation}. We obtain
\[
\begin{aligned}
R(x,y)&=2\pi^2\sum^{\infty}_{\ell=0}\sum^{\ell}_{m=-\ell}
Y_{\ell m}(\theta_1,\varphi_1)\overline{Y_{\ell m}(\theta_2,\varphi_2)}\\
&\quad\times\int^{\infty}_{\mathrm{i}}\frac{|\Gamma(\mathrm{i}\omega)|^2
P^{-1/2-\ell}_{-1/2+\mathrm{i}\omega}(\cosh(\sqrt{-K}\chi_1))
\overline{P^{-1/2-\ell}_{-1/2+\mathrm{i}\omega}(\cosh(\sqrt{-K}\chi_1))}}
{|\Gamma(\mathrm{i}\omega-\ell)|^2\sqrt{\sinh(\sqrt{-K}\chi_1)\sinh(\sqrt{-K}\chi_2)}}
\,\mathrm{d}\nu(\omega).
\end{aligned}
\]
Using the formula $\Gamma(z)=(z-1)\Gamma(z-1)$ and \eqref{eq:negativeR}, we
have
\[
f(\chi,\theta,\varphi)=2\sqrt{\pi}\sum^{\infty}_{\ell=0}\sum^{\ell}_{m=-\ell}
Y_{\ell m}(\theta,\varphi)\int^{\infty}_{\mathrm{i}\sqrt{-K}}
R_{k\ell}(\chi)\,\mathrm{d}Z_{\ell m}(k).
\]

The case of the flat model is similar, but this time
\[
V=\{\,([0,\infty),\ell,m)\colon\ell\in\mathbb{Z}_+,-\ell\leq m\leq\ell\,\}.
\]
The associated spherical functions for this case were known in implicit form
since \cite{Erdelyi1953b}. The group-theoretical interpretation has been
proposed by \cite{Vilenkin1957}. We have
\[
U^{\omega}_{\ell,0;0,0}(r)=\mathrm{i}^{\ell}\sqrt{2\ell+1}j_{\ell}(\omega r),
\]
and the spectral expansion has the form
\[
f(\chi,\theta,\varphi)=\pi\sqrt{2}\sum^{\infty}_{\ell=0}\sum^{\ell}_{m=-\ell}
Y_{\ell m}(\theta,\varphi)\int^{\infty}_0R_{k\ell}(\chi)\,\mathrm{d}Z_{\ell m}(k).
\]

Finally, in the case of the closed model
\[
V=\{\,(\omega,\ell,m)\colon\omega\in\mathbb{Z}_+,0\leq\ell\leq\omega,
-\ell\leq m\leq\ell\,\},
\]
$\mathfrak{V}$ is the $\sigma$-field of all subsets of $V$, and
$\mu(\omega,\ell,m)=\nu(\omega)$. The associated spherical functions are
known since \cite{Erdelyi1953b} and have the form
\[
U^{\omega}_{\ell,0;0,0}(r)=2^{\ell}\ell!\sqrt{\frac{(\omega-\ell)!(2\ell+1)}
{(\omega+\ell+1)!(\omega+1)}}\sin^{\ell}rC^{\ell+1}_{\omega-\ell}(\cos r),
\]
where $C^{\ell+1}_{\omega-\ell}(\cos r)$ are Gegenbauer polynomials. Using
the formula
\[
C^p_q(\cos r)=\frac{\sqrt{\pi}(p+2q-1)!\sin^{1/2-q}r}{2^{q-1/2}(q-1)!p!}
P^{1/2-q}_{p+q-1/2}(\cos r),
\]
the value $\omega=k/\sqrt{K}-1$, and \eqref{eq:positiveR}, we obtain
\[
U^{\omega}_{\ell,0;0,0}(r)=\sqrt{2\ell+1}R_{k\ell}(\chi)
\]
and
\[
f(\chi,\theta,\varphi)=2\sqrt{\pi}\sum^{\infty}_{\omega=0}\sum^{\omega}_{\ell=0}
R_{k\ell}(\chi)\sum^{\ell}_{m=-\ell}Y_{\ell m}(\theta,\varphi)Z_{k\ell m}.
\]

In what follows we suppose that the measure $\nu$ is absolutely continuous
with respect to the Plancherel measure. In other words,
\[
R(\chi)=
\begin{cases}
\frac{1}{\sqrt{-K}}\int^{\infty}_0\frac{\sqrt{-K}\sin(k\chi)}{k\sinh(\sqrt{-K}\chi)}
P(k)k^2\,\mathrm{d}k,&K>0,\\
\int^{\infty}_0\frac{\sin(k\chi)}{k\chi}P(k)k^2\,\mathrm{d}k,&K=0,\\
\frac{1}{K}\sum^{\infty}_{k/\sqrt{K}=1}\frac{\sqrt{K}\sin(k\chi)}{k\sin(\sqrt{K}\chi)}P(k)k^2,&K>0,
\end{cases}
\]
where $P(k)$ is the density of the measure $\nu$ with respect to the
Plancherel measure. $P(k)$ is called the \emph{spectral density} by
mathematicians and the \emph{power spectrum} by cosmologists. The power
spectrum $P(k)$ is the spherical Fourier transforms of the autocorrelation
function $R(\chi)$, while the latter is the inverse spherical Fourier
transform of the power spectrum. We have the following theorem.

\begin{theorem}\label{th:1}
The spectral expansion of the $\mathcal{G}$-invariant random field $f$ takes
the following form. In the case of the open model
\begin{equation}\label{eq:open}
f(\chi,\theta,\varphi)=2\sqrt{\pi}\sum^{\infty}_{\ell=0}\sum^{\ell}_{m=-\ell}
Y_{\ell m}(\theta,\varphi)\int^{\infty}_0
R_{k\ell}(\chi)k\sqrt{P(k)}\,\mathrm{d}W_{\ell m}(k).
\end{equation}
In the case of the flat model
\[
f(\chi,\theta,\varphi)=\pi\sqrt{2}\sum^{\infty}_{\ell=0}\sum^{\ell}_{m=-\ell}
Y_{\ell m}(\theta,\varphi)\int^{\infty}_0R_{k\ell}(\chi)k\sqrt{P(k)}
\,\mathrm{d}W_{\ell m}(k).
\]
Finally, in the case of the closed model
\[
f(\chi,\theta,\varphi)=2\sqrt{\pi}\sum^{\infty}_{\omega=0}\sum^{\omega}_{\ell=0}
R_{\omega\ell}(\chi)\sum^{\ell}_{m=-\ell}Y_{\ell m}(\theta,\varphi)\omega
\sqrt{P(\omega)}W_{\omega\ell m},
\]
where $W_{\ell m}$ is a sequence of uncorrelated identically distributed
complex-valued scattered random measures on $[0,\infty)$ with Lebesgue
measure as their common control measure, and where $W_{\omega\ell m}$ is a
sequence of uncorrelated identically distributed random variables with zero
mean and unit variance.
\end{theorem}

\section{Spin fields}\label{sec:spin}

In order to motivate the introducing of spin fields, consider the following
example. Let $\rho(t,x)$, $x=(\chi,\theta,\varphi)\in X$ be the density of
the non-relativistic matter. Currently, only the baryonic part of $\rho$ is
observable. Define the \emph{fractional overdensity} as
\[
\delta(t,x)=\frac{\rho(t,x)-\rho_{M,t}}{\rho_{M,t}}.
\]
The \emph{Newton's gravitational potential}, $\phi$, is the solution to
Poisson's equation
\[
\Delta_K\phi(t,x)=\frac{4\pi G}{c^2}\rho(t,x).
\]
Consider a source of electromagnetic waves at the comoving distance $\chi$
with angular coordinates $\theta$ and $\varphi$. The \emph{effective lensing
potential} of the source is defined as
\[
\psi(\chi,\theta,\varphi)=\frac{2}{f_K(\chi)}\int^{\chi}_0\frac{f_K(\chi-r)}{f_K(r)}
\phi(t_0-t_L,r,\theta,\varphi)\,\mathrm{d}r.
\]
This quantity is not observable. However, the \emph{convergence}
$\kappa=(1/2)\eth^*\eth\psi$ and the \emph{shear} $\gamma=(1/2)\eth^2\psi$
are observable now, while the \emph{first gravitational flexion}
$F=(1/2)\eth\eth^*\eth\psi$ and the \emph{second gravitational flexion}
$G=(1/2)\eth^3\psi$ may become observable in the near future.

The symbol $\eth$ denote the \emph{spin raising operator}. There exist
several equivalent definitions of this object. Consider the one connected to
representation theory.

Let $s\in\mathbb{Z}$ be an integer number. The map $s\mapsto(\varphi\mapsto
e^{\mathrm{i}s\varphi})$, $\varphi\in\mathrm{SO}(2)$ determines a one-to-one
correspondence between $\mathbb{Z}$ and the set of the irreducible unitary
representations of the group $\mathcal{H}=\mathrm{SO}(2)$. Consider the
Cartesian product $\mathrm{SO}(3)\times\mathbb{C}$ and introduce the
following equivalence relation: the points $(g_1,z_1)$ and $(g_2,z_2)$ are
equivalent if there exists an element $h=h_{\varphi}\in\mathcal{H}$ such that
$g_2=g_1h$ and $z_2=e^{-\mathrm{i}s\varphi}z_1$. Call the quotient space
$E_s$. The projection map $\pi(g,z)=g\mathcal{H}$ maps $E_s$ onto the
homogeneous space $\mathrm{SO}(3)/\mathrm{SO}(2)$, which is the
two-di\-men\-si\-o\-nal sphere, $\mathbb{S}^2$. The inverse image of each
point in $\mathbb{S}^2$ is a copy of $\mathbb{C}$. The triple
$(E_s,\pi,\mathbb{S}^2)$ is a homogeneous line bundle over $\mathbb{S}^2$.

Let $f\colon\mathbb{S}^2\to E_s$ be a \emph{section} of the above bundle,
i.e., $\pi\circ f$ is the identity map on $\mathbb{S}^2$. The group
$\mathrm{SO}(3)$ acts on $E_s$ by $g(g_0,z)=(gg_0,z)$. This action identifies
the fibers over any two points of the base space. Therefore, we can define
the square-integrable sections by the condition
\[
\int^{2\pi}_0\int^{\pi}_0|f(\theta,\varphi)|^2\sin\theta\,\mathrm{d}\theta
\,\mathrm{d}\varphi<\infty.
\]
The elements of this space are called the square-integrable functions of
\emph{spin} $s$.

The \emph{induced representation} $U$ of the group $\mathrm{SO}(3)$ acts in
the space $L_2(E_s)$ of the square-integrable sections by
\[
U(g)f(n)=f(g^{-1}n),\qquad g\in\mathrm{SO}(3),\quad n\in\mathbb{S}^2.
\]
The \emph{Frobenius reciprocity} states that the multiplicity of the
irreducible representation $U_{\ell}$ in $U$ is equal to the multiplicity of
the representation $s$ in $U_{\ell}$. The latter is equal to $1$ if
$|s|\leq\ell$ and to $0$ otherwise. It follows that the induced
representation $U$ is equivalent to the direct sum of irreducible components
$U_{\ell}$, $\ell=|s|$, $|s|+1$, \dots.

The basis in the space $L_2(E_0)$ consists of spherical harmonics
\eqref{eq:spherical} (of spin 0). The basis in $L_2(E_s)$ consists of
\emph{spherical harmonics of spin} $s$:
\[
{}_sY_{\ell m}(\theta,\varphi)=\sqrt{\frac{2\ell+1}{4\pi}}D^{(\ell)}_{-m,-s}(\varphi,\theta,0).
\]
The last display and \eqref{eq:Wigner} imply the following formula:
\[
\begin{aligned}
{}_sY_{\ell m}(\theta,\varphi)&=e^{\mathrm{i}m\varphi}(-1)^m
\sqrt{\frac{(2\ell+1)(\ell+m)!(\ell-m)!}{4\pi(\ell+s)!(\ell-s)!}}
\sin^{2\ell}(\theta/2)\\
&\quad\times\sum^{\min\{\ell+m,\ell-s\}}_{u=\max\{0,m-s\}}\binom{\ell-s}{u}
\binom{\ell+s}{u-m+s}(-1)^{\ell-u-s}\cot^{2u-m+s}(\theta/2).
\end{aligned}
\]

Consider the operator $\eth$ (in fact, a family of operators) acting onto a
spherical harmonic of spin $s$ by
\[
\eth{}_sY_{\ell m}=\sqrt{(\ell-s)(\ell+s+1)}{}_{s+1}Y_{\ell m}.
\]
This action may be extended by linearity to the linear subset
$\mathcal{D}_{\eth}$ of the space $L_2(E_s)$ that consists of the functions
${}_sf$, for which their Fourier coefficients
\[
{}_sf_{\ell m}=\int_{\mathbb{S}^2}{}_sf(\theta,\varphi)
\overline{{}_sY_{\ell m}(\theta,\varphi)}\sin\theta\,\mathrm{d}\theta
\,\mathrm{d}\varphi
\]
satisfy
\[
\sum^{\infty}_{\ell=|s|+1}(\ell-s)(\ell+s+1)|{}_sf_{\ell m}|^2<\infty.
\]
It is easy to check, that
\[
\eth=s\cot\theta-\frac{\partial}{\partial\theta}-\frac{\mathrm{i}}{\sin\theta}
\frac{\partial}{\partial\varphi}.
\]

The conjugate operator, $\eth^*$, acts onto a spherical harmonic of spin $s$
by
\[
\eth^*{}_sY_{\ell m}=-\sqrt{(\ell+s)(\ell-s+1)}{}_{s-1}Y_{\ell m}.
\]
By obvious reasons, $\eth$ is called the \emph{spin raising operator}, while
$\eth^*$ is called the \emph{spin lowering operator}. Moreover, the
restrictions of these operators to the space $H_{\ell}$, where the
irreducible representation $U_{\ell}$ acts, are intertwining operators
between equivalent representations.

Assume that the effective lensing potential of the sources of radiation is a
single realisation of a random field. By the cosmological principle, the
above field must be $\mathcal{G}$-invariant. In particular, its spectral
expansion is determined by Theorem~\ref{th:1}. Consider the case of the
convergence field in an open model. Assume that the power spectrum $P(k)$
satisfies the following condition:
\[
\sum^{\infty}_{\ell=1}\ell(\ell+1)(2\ell+1)\int^{\infty}_0|R_{k\ell}(\chi)|^2
k^2P(k)\,\mathrm{d}k<\infty.
\]
Apply operator $(1/2)\eth^*\eth$ to the expansion \eqref{eq:open} term by
term:
\[
\kappa(\chi,\theta,\varphi)=-\sqrt{\pi}\sum^{\infty}_{\ell=1}\sum^{\ell}_{m=-\ell}
\ell(\ell+1)Y_{\ell m}(\theta,\varphi)\int^{\infty}_0R_{k\ell}(\chi)k\sqrt{P(k)}
\,\mathrm{d}W_{\ell m}(k).
\]
The series in the right hand side converges in mean square and determines the
spectral expansion of a random field. Note that the above expansion cannot be
written in the form of \eqref{eq:open}. Therefore, the random field
$\kappa(\chi,\theta,\varphi)$ is not $\mathrm{SO}_0(1,3)$-invariant. However,
the restriction of the field $\kappa$ to the sphere of a fixed comoving
radius $\chi_0$ centred at the fundamental observer is an isotropic random
field on the sphere. In other words, the above restriction is an isotropic
random section of the line bundle $(E_0,\pi,\mathbb{S}^2)$.

The heuristic arguments of the previous paragraph suggest the following. The
restriction of the convergence field to the sphere of a fixed comoving radius
$\chi_0$ centred at the fundamental observer is an isotropic random section
of the line bundle $(E_0,\pi,\mathbb{S}^2)$, that of the first gravitational
flexion field is an isotropic random section of the line bundle
$(E_1,\pi,\mathbb{S}^2)$, that of the shear field is an isotropic random
section of the line bundle $(E_2,\pi,\mathbb{S}^2)$. Finally, the restriction
of the second gravitational flexion field to the sphere of a fixed comoving
radius $\chi_0$ centred at the fundamental observer is an isotropic random
section of the line bundle $(E_3,\pi,\mathbb{S}^2)$. This means that the
expected value $\mathsf{E}[X(\chi_0,\mathbf{n})]$,
$X\in\{\kappa,\gamma,F,G\}$, $\mathbf{n}\in S^2$, is constant, while the
covariance function
$\mathbf{E}[X(\chi_0,\mathbf{n}_1)\overline{X(\chi_0,\mathbf{n}_2)}]$
satisfies the following condition:
\[
\mathbf{E}[X(\chi_0,g\mathbf{n}_1)\overline{X(\chi_0,g\mathbf{n}_2)}]
=\mathbf{E}[X(\chi_0,\mathbf{n}_1)\overline{X(\chi_0,\mathbf{n}_2)}],\qquad
g\in\mathrm{SO}(3).
\]

Spectral expansion of an isotropic random field in the trivial bundle
$(E_0,\pi,\mathbb{S}^2)$ goes back to \citet{obukhov}, that in the bundle
$(E_2,\pi,\mathbb{S}^2)$ is known to cosmologists since \citet{zaldarriaga}.
The rigourous mathematical theory of isotropic random fields in vector
bundles was proposed by \citet{Geller2010} and \citet{Malyarenko2011a}. In
particular, the spectral expansion of an isotropic random section
$X(\theta,\varphi)$ of the line bundle $(E_s,\pi,\mathbb{S}^2)$,
$s\in\mathbb{Z}$, has the form
\[
X(\theta,\varphi)=\sum^{\infty}_{\ell=s}\sum^{\ell}_{m=-\ell}a_{s,\ell m}\,
{}_sY_{\ell m}(\theta,\varphi)
\]
where $\mathsf{E}[a_{s,\ell m}]=0$ unless $s=\ell=0$ and
\[
\begin{aligned}
\mathsf{E}[a_{s,\ell m}\overline{a_{s,\ell'm'}}]&=\delta_{\ell\ell'}\delta_{mm'}C_{s\ell},\\
\sum^{\infty}_{\ell=s}(2\ell+1)C_{s\ell}&<\infty.
\end{aligned}
\]

It follows that there exists a sequence of uncorrelated mean square
continuous stochastic processes $a_{s,\ell m}(\chi)$ such that
\begin{equation}\label{eq:spin}
X(\chi,\theta,\varphi)=\sum^{\infty}_{\ell=s}\sum^{\ell}_{m=-\ell}a_{s,\ell m}
(\chi){}_sY_{\ell m}(\theta,\varphi).
\end{equation}
Moreover, $\mathsf{E}[a_{s,\ell m}(\chi)]=0$ unless $s=\ell=m=0$ and
\begin{equation}\label{eq:expansion1}
\mathsf{E}[a_{s,\ell m}(\chi_1)\overline{a_{s,\ell'm'}(\chi_2)}]
=\delta_{\ell\ell'}\delta_{mm'}C_{s\ell}(\chi_1,\chi_2)
\end{equation}
with
\begin{equation}\label{eq:expansion2}
\sum^{\infty}_{\ell=s}(2\ell+1)C_{s\ell}(\chi,\chi)<\infty.
\end{equation}
As usual, the comoving distance $\chi$ runs from $0$ to $\pi/\sqrt{K}$ when
$K>0$, and from $0$ to $\infty$ otherwise.

The correlation function of the random field \eqref{eq:spin} is
\[
\begin{aligned}
R(\chi_1,\theta_1,\varphi_1,\chi_2,\theta_2,\varphi_2)
&=\mathsf{E}[X(\chi_1,\theta_1,\varphi_1)\overline{X(\chi_2,\theta_2,\varphi_2)}]\\
&=\sum^{\infty}_{\ell=s}C_{s\ell}(\chi_1,\chi_2)\sum^{\ell}_{m=-\ell}
{}_sY_{\ell m}(\theta_1,\varphi_1)\overline{{}_sY_{\ell
m}(\theta_2,\varphi_2)}.
\end{aligned}
\]
To simplify this formula, use the addition theorem for spin weighted
spherical harmonics fromz \citet[equation~(A4.57)]{durrer}:
\[
\sqrt{\frac{4\pi}{2\ell+1}}\sum^{\ell}_{m'=-\ell}{}_sY_{\ell m'}(\theta_2,\varphi_2)
\overline{{}_{-m}Y_{\ell m'}(\theta_1,\varphi_1)}={}_sY_{\ell m}(\beta,\alpha)
e^{-\mathrm{i}s\gamma},
\]
where $(\alpha,\beta,\gamma)$ are the Euler angles of the rotation which
first rotates $(\theta_2,\varphi_2)$ into the north pole and then north pole
into $(\theta_1,\varphi_1)$. In fact, $\beta$ is the angle between
$(\theta_1,\varphi_1)$ and $(\theta_2,\varphi_2)$, $\alpha$ is the angle
between the great arc connecting the above two points and the meridian
passing through $(\theta_1,\varphi_1)$, and $\gamma$ is the angle between the
great arc connecting the above two points and the meridian passing through
$(\theta_2,\varphi_2)$, see Fig.~1 in \citet{Ng1999}. Substitute the value
$m=-s$, use formula
\begin{equation}\label{eq:phase}
{}_sY_{\ell m}(\beta,\alpha)=e^{-\mathrm{i}m\alpha}{}_sY_{\ell m}(\beta,0),
\end{equation}
and interchange indices $1$ and $2$. We obtain
\[
\sum^{\ell}_{m=-\ell}{}_sY_{\ell m}(\theta_1,\varphi_1) \overline{{}_sY_{\ell
m}(\theta_2,\varphi_2)}=\sqrt{\frac{2\ell+1}{4\pi}}{}_sY_{\ell
-s}(\beta,0)e^{-\mathrm{i}s(\alpha+\gamma)},
\]
and
\begin{equation}\label{eq:corr}
R(\chi_1,\theta_1,\varphi_1,\chi_2,\theta_2,\varphi_2)
=\frac{1}{2\sqrt{\pi}}\sum^{\infty}_{\ell=s}C_{s\ell}(\chi_1,\chi_2)\sqrt{2\ell+1}
{}_sY_{\ell -s}(\beta,0)e^{-\mathrm{i}s(\alpha+\gamma)}.
\end{equation}

To calculate $C_{s\ell}(\chi_1,\chi_2)$, put $\theta_1=\varphi_2=0$ and
$\theta_2=\beta$ in the last display. In other words, let
$(\theta_1,\varphi_1)$ be the north pole. Then $\alpha=\gamma=0$ and
\[
R(\chi_1,0,0,\chi_2,\beta,0)=\frac{1}{2\sqrt{\pi}}\sum^{\infty}_{\ell=s}
C_{s\ell}(\chi_1,\chi_2)\sqrt{2\ell+1}{}_sY_{\ell -s}(\beta,0).
\]
The spin weighted spherical harmonics ${}_sY_{\ell m}(\theta,\varphi$,
$\ell\geq|s|$, $-\ell\leq m\leq\ell$ form an orthonormal basis in the space
of square integrable functions on the sphere $\mathbb{S}^2$. It follows from
\eqref{eq:phase} that the functions $\sqrt{2\pi}{}_sY_{\ell -s}(\theta,0)$
form an orthonormal basis in the space $L^2([0,\pi],\sin\theta\,d\theta)$.
Multiply both hand sides of the last display by $\overline{{}_sY_{\ell
-s}(\beta,0)}$ and integrate over $[0,\pi]$ with respect to the measure
$\sin\beta\,\mathrm{d}\beta$. We obtain
\begin{equation}\label{eq:C}
C_{s\ell}(\chi_1,\chi_2)=\frac{4\pi^{3/2}}{\sqrt{2\ell+1}}\int^{\pi}_0
R(\chi_1,0,0,\chi_2,\beta,0)\overline{{}_sY_{\ell -s}(\beta,0)}\sin\beta\,\mathrm{d}\beta.
\end{equation}

Multiply both hand sides of equation \eqref{eq:spin} by
$\overline{{}_sY_{\ell m}(\theta,\varphi)}$ and integrate over
$\mathbb{S}^2$. We obtain
\[
a_{s,\ell m}(\chi)=\int_{\mathbb{S}^2}X(\chi,\theta,\varphi)
\overline{{}_sY_{\ell m}(\theta,\varphi)}\sin\theta\,\mathrm{d}\theta
\,\mathrm{d}\varphi.
\]
It follows from mean square continuity that
\[
\lim_{\chi\downarrow 0}\mathsf{E}[|a_{s,\ell
m}(\chi)|^2]=0,\qquad(s,\ell,m)\neq(0,0,0),
\]
therefore
\begin{equation}\label{eq:aslm}
C_{s\ell}(0,\chi)=0,\qquad\ell\neq 0.
\end{equation}

We arrive at the following theorem.

\begin{theorem}\label{th:2}
The spin cosmological fields have the form \eqref{eq:spin}, where $s=0$ for
the convergence, $s=1$ for the first gravitational flexion, $s=2$ for the
shear, and $s=3$ for the second gravitational flexion. The sequence of
stochastic processes $a_{s,\ell m} (\chi)$ satisfies \eqref{eq:expansion1}
and \eqref{eq:expansion2}. The correlation function of a spin cosmological
field has the form \eqref{eq:corr}, where $C_{s\ell}(\chi_1,\chi_2)$ is a
sequence of nonnegative-definite kernels satisfying \eqref{eq:expansion2} and
\eqref{eq:aslm}. The kernels $C_{s\ell}(\chi_1,\chi_2)$ may be calculated by
\eqref{eq:C}.
\end{theorem}

Theorem~\ref{th:2} in the case of $s=0$ is due to \citet{Yadrenko1983}.

\bibliographystyle{abbrvnat}
\bibliography{cosmo}

\end{document}